\documentclass{amsart}
\usepackage{amsmath,amssymb,amsthm}
\usepackage[rokicki]{boxedeps}
\usepackage{graphicx}
\HideDisplacementBoxes

\newtheorem{thm}{Theorem}[section]

\newtheorem{cor}{Corollary}[section]

\numberwithin{equation}{section}
\def\A{\Bbb A}
\def\Z{\Bbb Z}

\def\R{\Bbb R}

\def\P{\Bbb P}
\def\C{\Bbb C}

\def\flabel#1{\ifmmode #1\else$ #1$\fi}

\let\angle\undefined
\newsymbol\angle 115C
\def\degrees{\ifmmode^\circ\else$^\circ$\fi}

\def\pict\includegraphics

\def\picture\includegraphics

\catcode`\@=11
\def\begintable{\bigbreak
\begingroup\interlinepenalty=1000
\def\bline{\noalign{\vskip\tableskip\nobreak}}}

\begin{document}

\title[Curves in cages]{curves in cages: 
 an algebro-geometric zoo}

\author{Gabriel Katz}

\address{Brandeis University, Waltham, MA 02454-9110}

\email{\hfil\break gabrielkatz@rcn.com} 

\maketitle

\section{Introduction.}

An algebraic plane curve is the solution set of a polynomial equation $P(x, y) = 0$, where $x$ and $y$
are real or complex variables. By definition, the {\it degree} of the curve is the degree of the polynomial
$P(x, y)$.\footnote {According to this  point of view adapted by algebraic geometers, the degree is associated with the equation rather than with its solution set: for example, the equation $(x^2 + y^2 - 1)^2 = 0$ defines a ``double circle'" of degree four.}  When $P(x, y)$ is a product of two non-constant polynomials over a given number field,  the curve is called {\it reducible}; otherwise, it is  {\it irreducible}. \smallskip

This paper is concerned with families of plane algebraic curves that contain
a given and quite  special  set of points $\mathcal X$.  We  focus on the case in which 
the set $\mathcal X$ is formed by transversally intersecting pairs of 
lines selected from two given finite families.   
The union of all lines from both families is called a 
\emph{cage} (this notion of  cage will be made more precise later), and the 
intersection $\mathcal X$ consists of points at which a line from the first family 
intersects a line from the second. 
The points of $\mathcal X$ are called the {\it nodes} of the cage.

This  is a particular case of a more general problem.
Let $\mathcal X$ be the intersection set of two plane algebraic curves 
$\mathcal D$ and $\mathcal E$ that do not share a common component, that is, do not contain a common irreducible curve. If $d$ and $e$  denote the degrees of $\mathcal D$ and $\mathcal E$, respectively, then  $\mathcal X$  consists  of at most $d \cdot e$ points. When the cardinality 
of $\mathcal X$ is exactly $d \cdot e$, $\mathcal X$ is called a 
\emph{complete intersection} (complete intersections have many nice 
properties). 
How does one describe polynomials of degree at most $k$ that vanish on a complete 
intersection $\mathcal X$ or on its 
subsets? This problem has a glorious history (see [{\bf 1}], [{\bf 2}], or [{\bf 3}]) and its 
generalizations are a subject of active and exciting research 
[{\bf 4}], [{\bf 5}], [{\bf 6}], [{\bf 7}], [{\bf 8}], [{\bf 9}].  
When  $k$ is much larger then $d\cdot e$, points of $\mathcal X$ impose \emph{independent} restrictions on polynomials of degree $k$.  However, for small $k$ these restrictions
fail to be independent. In such cases, we can look for  maximal 
subsets of $\mathcal X$ that impose independent constraints.

In this article, we deal with the special case of this classical problem suggested earlier, namely,
the case in which both plane curves $\mathcal D$ and $\mathcal E$ are simply 
\emph{unions of lines} and the union $\mathcal D \cup \mathcal E$ is the $(d \times e)$-cage 
in question. Note that $\mathcal D$ is the zero set of a product of $d$ linear polynomials, while $\mathcal E$ is the zero set of a product of $e$ linear polynomials. Hence, the degrees of 
$\mathcal D$ and $\mathcal E$ are $d$ and $e$, respectively.
In order to simplify our terminology, we color the lines from $\mathcal D$ 
\emph{red} and the ones from $\mathcal E$ \emph{blue}.
\smallskip

The case of cages is amenable to elementary methods that presume only a 
modest familiarity with algebraic geometry. Moreover, 
its beautiful geometric applications are the true focus 
of our exposition. These applications can be viewed as natural 
generalizations of classical theorems in the foundations of projective 
geometry. The reader familiar with Pascal's theorem can get a feel for 
the nature of these generalizations by 
glancing at Figures 6 and 7 and comparing them with 
the classical Pascal diagram in Figure 1. 
\smallskip

The zoo of algebraic curves attached to the nodes of a cage is a 
microcosm of classical algebraic geometry. We invite the Monthly readership to join us on a tour of its feature attractions. 
\bigskip

My interest in the subject was sparked by N. B. Vasiliev's 
engaging paper ``Pascal's Hexagrams and Cubic Curves" [{\bf 13}].
The striking and well-known connection
between classical theorems of projective geometry and
the algebraic geometry of cubic plane curves
 was a revelation to me. Equally striking was the
elegance of the argument that established
this connection. This article was written in the afterglow of this 
private epiphany. \bigskip

\section{Pascal's Mystic Hexagram.}

In 1640, a the sixteen-year-old Blaise Pascal discovered a
remarkable property of a hexagon inscribed in a
 circle. In a diagram generated by the hexagon, 
three specific intersection points always happen to be collinear!
He called the configuration `` The Mystic Hexagram," made
fifty posters of it, and mailed them to  fellow scientists. 

Shortly thereafter, Pascal realized that a similar observation holds for a hexagon inscribed 
in an ellipse. In fact, the amazing collinearity is preserved under central and parallel 
projections of his diagram (see Figure 1): a projection maps lines to lines, the circle 
is transformed into a quadratic curve,  and the inscribed hexagon is mapped onto a 
hexagon inscribed in that curve.
 
Pascal's mystic hexagram was a fundamental result in geometry unknown
to the classical Greek school. In the words of Fermat, ``we
learned that the ancient Greeks
did not know everything about geometry."
Along with Desargues's theorem
(discovered four years earlier), Pascal's theorem (see [{\bf 11}, Corollary  3.15]) gave birth to a new branch of nonmetrical geometry that we now call  \emph{projective geometry}.

\begin{figure}[ht]
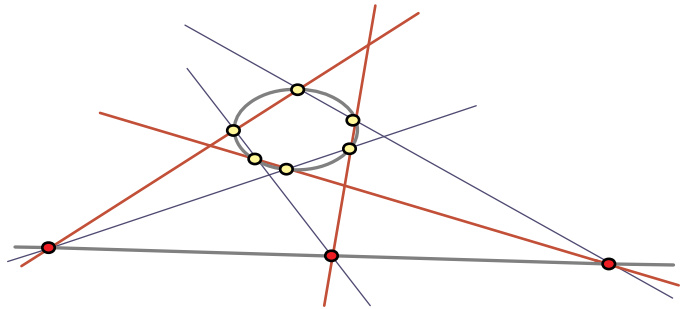

\centerline{\BoxedEPSF{Pascal scaled 500}}
\bigskip
\caption{Pascal's mystic hexagram.}
\end{figure} 

\begin{thm} {\bf (Pascal's Mystic Hexagram)}. 
Given a hexagon inscribed in a
quadratic curve $\mathcal Q$, color its edges alternately blue
and red and extend them.
Generically, the resulting three red lines
intersect the three blue lines at three 
points different from the six vertices of the hexagon.
These three points must be \emph{collinear}. 
\end{thm}

The three blue and three red lines generated by the
hexagon form  a $(3 \times 3)$-cage. In general, two triples of
lines  produce nine nodes (intersection points), 
each of which belongs to a unique pair of one red line and one blue line.

In the mystic hexagram, six of the nine nodes lie on the quadratic
curve $\mathcal Q$,  while the remaining nodes (by the Pascal theorem)
must lie on a line $\mathcal L$. Consider the union 
$\mathcal C = \mathcal Q \cup \mathcal L$. The set $\mathcal C$ is an
example of a reducible \emph{cubic} curve---it is the
zero set of the product of two polynomials in two variables, one quadratic and
one linear.

We can now formulate a well-known corollary of Pascal's theorem:

\begin{thm}
All nine nodes of the  $(3\times 3)$-cage
generated by a hexagon inscribed in a quadratic curve
$\mathcal Q$ lie on a reducible cubic curve 
$\mathcal C = \mathcal Q \cup \mathcal L$. \qed
\end{thm}

Of course, the line $\mathcal L$ is determined by any pair $a$ and $b$ 
from the three new nodes $a$, $b$, and $c$ of the cage.
Thus, eight nodes of the cage---the six vertices of the
inscribed hexagon and the new pair $\{a, b\}$---always
belong to some cubic curve $\mathcal C$.  What happens
to the remaining \emph{ninth} node $c$? According to the
Pascal theorem, the ninth node must lie on $\mathcal C$ as well!
This reformulation of Pascal's theorem leads one to
wonder whether the theorem might be a general cubic phenomenon.
Might it  hold for a hexagon inscribed in a general cubic curve? Perhaps the
reducibility of $\mathcal C = \mathcal Q \cup \mathcal L$ 
is merely an accident.
The following classical proposition (see  [{\bf 12}, Corollary 2.7]) validates these insights (Figure 2): 
\begin{thm}{\bf (Cage Theorem for Cubics).}
Any cubic curve $\mathcal C$ that passes through eight nodes of a
$(3\times 3)$-cage must pass through the ninth node.
\end{thm}
\begin{figure}[ht]
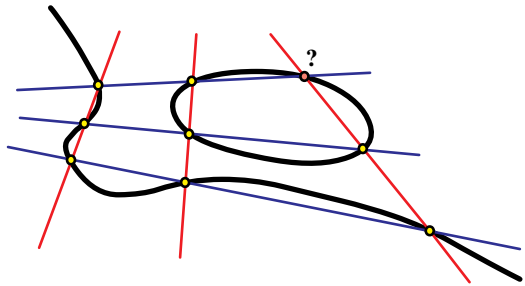

\centerline{\BoxedEPSF{3CageTh scaled 500}}
\bigskip
\caption{The cage theorem for cubics.} 
\end{figure}
A natural generalization of Theorem 2.3 is a theorem of Chasles [{\bf 3}].  It claims
that if $\mathcal X$ is a complete intersection of 
two cubic curves, then any cubic curve that contains eight points from $\mathcal X$ 
will include its ninth point as well.\bigskip 

Before proceeding, we issue one caveat to the reader. Although the diagrams in this paper suggest the case of real algebraic curves, all our 
algebraic arguments and the results make sense---and  are even simpler---for \emph{complex} 
algebraic curves (defined by complex polynomials). 
\bigskip

At first glance, Theorem 2.3 appears to describe an esoteric fact.
However, it reflects a wonderful intrinsic structure
shared by nonsingular cubic curves (called \emph{elliptic curves}). 
It turns out that any elliptic $\mathcal C$
curve has the structure of an \emph{Abelian group} (see [{\bf 11}, 3.14]). 

We recall for the reader how the group operation $(x, y)  \longmapsto x + y$ is defined for an elliptic curve $\mathcal C$. Here we view  $\mathcal C$ as residing in a projective plane (see section 3 for a short discussion of projective curves and spaces). If $x$ and $y$ are distinct points of $\mathcal C$, then
there is a unique point $z$ of $\mathcal C$ with the property that $x, y$, and $z$ are collinear;  if
$x = y$, then the tangent line to $\mathcal C$ at $x$ hits the curve at a single point $z$.
Next, we designate one point on $\mathcal C$ as the additive identity element $e$. 
 Once $e$ is chosen, $z$ must play the role of $-(x + y)$. In other words,  $x + y$ is uniquely  
determined by the property:  $x + y, z$, and $e$ are collinear points. Evidently, for each $x$ in 
$\mathcal C$ we have $x + e = x$. When there is one group structure  on a set, there are of course many others obtained from the given structure by conjugation with a translation. A translation sends the identity element $e$ to another element, which becomes the identity element of the new structure. It is customary to pick one of the inflection points of $\mathcal C$ for the role of $e$.

In the language of the group structure on an elliptic curve, the cage theorem becomes a statement
about the \emph{associativity} of the binary group operation! It's a  subtle interpretation.
Figure 3  indicates how it works.\bigskip
\begin{figure}[ht]
\centerline{\BoxedEPSF{associativity scaled 600}}
\bigskip
\caption{ }
\end{figure}

It is striking to realize that every classical theorem
of planar projective geometry is, in a sense, an
implication of  Theorem 2.3 for cubics, which in
turn reflects the associativity of point addition on
 nonsingular cubic curves! So elliptic curves and planar
projective geometry are intimately related.

One might wonder whether a different, more direct
generalization of Pascal's theorem is valid:
\emph{Is it true that any alternately-colored hexagon inscribed in a cubic
curve produces a $(3\times 3)$-cage all of whose nodes lie on the curve?}
A simple count of dimensions rules out this possibility.
For any cubic curve $\mathcal C$ the variety of 
all  hexagons inscribed in $\mathcal C$ is six-dimensional. One can show that 
the family of hexagons whose $(3\times 3)$-cages have all their nodes on
$\mathcal C$ constitutes a four-dimensional variety.\footnote{By Theorem 2.3, 
a generic quadrilateral inscribed in $\mathcal C$ gives rise to a unique $(3 \times 3)$-cage 
whose nodes belong to $\mathcal C$.} Thus, for an arbitrary $\mathcal C$, not
any inscribed hexagon will do---it is important for $\mathcal C$ to have a \emph{quadratic} 
component $\mathcal Q$ in which  the hexagons  can be inscribed! In particular, not any hexagon inscribed in a union of a quadratic curve with a line will generate a cage that belongs to that union.

%%%%%%%%%%%%%%%%%%%%%%%%%

\section{high degree curves in cages.}

We aim in this article to generalize some of the observations made in the previous section for
cubic curves to curves of any degree $d$ that pass through the nodes of a
$(d \times e)$-cage, where $e \leq d$. Most of our results can be 
derived from the Bacharach duality theorem [{\bf 1}], [{\bf 5}]  (see also Theorem 4.1, 
as outlined in section 4). 
However, our goal here is to replace the powerful machinery of algebraic geometry 
with more elementary considerations that rely only on the divisibility of polynomials.  \smallskip

Our arguments are identical 
for real or complex curves, meaning for curves over a ground field $\A$ with either $\A = \R$ or   
$\A = \C$. They apply both to curves in the affine $xy$-plane $\A^2$ and to curves in the projective plane $\P^2$ with homogeneous coordinates $[x : y : z]$. Recall that the points of $\P^2$ are the 
\emph{proportionality classes} of triples $(x, y, z)$ different from $(0, 0, 0)$ 
(i.e., the equivalence classes $[x : y : z]$ of $\A^3 \setminus \{(0, 0, 0)\}$ modulo the relation 
$``\sim": \;(x, y, z) \sim (\lambda x, \lambda y, \lambda z)$ for any 
$(x, y, z)$ in $\A^3 \setminus \{(0, 0, 0)\}$ and any $\lambda$ in $\A^\ast$, the multiplicative group of nonzero elements of $\A$). 
Depending on the context, the numbers $x, y$, and  $z$ can be real or complex. 
For the most part, our notations are neutral with respect to the choice 
of  the ground field $\A$.\smallskip

The introduction of the projective plane simplifies the intersection theory 
of algebraic curves. For example, in contrast with the plane $\A^2$, 
in the projective plane $\P^2$ any two lines have an intersection point.
\smallskip 

To put matters into context, we  recall for readers a few basic notions. 
Each polynomial  $P(x, y)$ of degree at most $d$ gives 
rise to a unique homogeneous polynomial in $x, y$, and $z$ of degree $d$. 
This is done simply by replacing each monomial $x^a y^b$ in $P(x, y)$ with the monomial 
$x^a y^b z^{d - a -b}$.

An \emph{affine (algebraic) curve} in $\A^2$ is the zero set of a single
polynomial in the variables $x$ and $y$, while
a \emph{projective curve} in the projective plane
$\P^2$ is the zero set of a \emph{homogeneous}
polynomial in $x, y$, and $z$. The degree of a curve is the degree
of its defining polynomial. Two polynomials that are proportional
define the same curve. Therefore, we interpret the set
of proportionality classes of polynomials in $x$ and $y$ of
degree $d$ as the set of plane affine curves of degree $d$.  Similarly, 
the set of proportionality classes of homogeneous polynomials in $x, y$, and $z$ of
degree $d$ may be seen as the set of projective curves of degree $d$ in $\P^2$.
For curves over the real numbers this is a
simplistic model (many polynomials have empty zero sets), 
but over the complex numbers it is more than adequate.
\smallskip

A polynomial $P(x, y)$ of degree not larger than $d$ has 
$$1 + 2 + 3 +\; \cdots \; + d = \hfil\break (d + 1)(d + 2)/2$$ coefficients.
Such polynomials, taken up to proportionality, constitute an $n$-dimensional projective
space  $\P_\ast^n$, where $$n =  [(d + 1)(d + 2)/2] - 1 = (d^2 + 3d)/2.$$  
Its points are the proportionality classes of sequences formed by the coefficients of $P(x, y)$.

When $k < d$, the polynomials in $x$ and $y$ of degree at most $k$ are contained in  the 
set of polynomials of degree at most $d$. This gives rise to a nested family of 
projective subspaces $\{\P_\ast^{(k^2 + 3k)/2}\}_{0 \leq k \leq d}$. 
For example, the family
$$ \P_\ast^0 \subset \P_\ast^2 \subset \P_\ast^5 \subset
\P_\ast^9 \subset \P_\ast^{14} \subset \P_\ast^{20}$$
corresponds to polynomials $P(x, y)$ of degrees 
$0$ and at most $1,  2,  3,  4,$ and $5$, respectively. 
The affine curves of degree $d$ in $\A^2$ form an open set
$$\P_\ast^{(d^2 + 3d)/2}\; \setminus \; \P_\ast^{[(d - 1)^2 + 3(d - 1)]/2}$$
in the space $\P_\ast^{(d^2 + 3d)/2}$. 
For instance, in this model the cubic curves in $\P^2$ are represented by points of
$\P_\ast^9$, while the cubic curves in $\A^2$ are represented by
the set  $\P_\ast^9 \setminus \P_\ast^5$.

As  $d$ increases, the members of  the universe of  plane curves of degree $d$ that contain  all the nodes of a $(d \times d)$-cage become more and more rare.
Indeed, the dimension of the variety
of plane curves of degree $d$ is a quadratic function of
$d$, whereas  the dimension of the subvariety of  such 
curves that contain the nodes of a $(d \times d)$-cage grows only linearly in $d$. 
Still, among the caged curves there are a few interesting beasts. 
Consider, for instance, the Fermat curve $\mathcal{F}$ given 
in homogeneous coordinates $[x : y : z]$ by the equation 
$x^d + y^d = z^d$, or in Cartesian (i.e. affine) coordinates  by 
the equation  $x^d + y^d = 1$. 
We notice that, over the complex numbers, $\mathcal{F}$ is a curve 
that contains the nodes of a $(d \times d)$-cage. 
Simply write the curve's equation in the form
  $(x^d - 1/2) +  (y^d - 1/2) = 0$ or in the form 
$ \prod_{\xi }\, (x - \xi ) + \prod_{\xi' }\, (y - \xi' ) = 0$, 
where $\xi$ and $ \xi'$ run over the $d$th roots of $1/2$. 
Thus, all the root pairs $(\xi , \xi')$ form the nodes of the complex 
\emph{Fermat cage}. 

Our next goal is to describe polynomials $P$ of degree $d$
that vanish at the nodes of a given $(d \times e)$-cage $K$, where $d \geq e$.
By forcing a polynomial $P$ to vanish at a
given point of a projective plane,
we are imposing a single linear restriction on its coefficients.
More points  produce further restrictions, and could
lower the dimension of the vector space formed by such polynomials. 
If the degrees of polynomials are not bounded from above, 
distinct points will impose independent constraints. However, 
when the degrees are bounded, new points often add only 
redundant constraints.  

By  way of example, one expects that requiring a polynomial $P(x, y)$ of degree
less then or equal to $d$ to vanish at $(d^2 + 3d)/2$ generic points should
nail it down (up to proportionality).
Since, a $(d \times d)$-cage has $d^2$ nodes,  some of them must impose
redundant restrictions on the coefficients of $P(x, y)$ (provided that $d \geq 3$).

The issue is: Which nodes of a $(d \times e)$-cage are redundant, and how can one
describe the variety of degree $d$ polynomials vanishing at the
nodes?
A few combinatorial definitions, linked to the notion of $k$-configurations in [{\bf 6}], 
will help us to depict the set of redundant nodes.  
As before, we mark the two sets of lines forming a cage with two colors.
Given any ordering of the red lines
$\{\mathcal R_1 , \mathcal R_2 , ... , \mathcal R_d \}$ and the blue 
lines $\{\mathcal B_1 , \mathcal B_2 , ... , \mathcal B_e \}$, we   
denote by $p_{ij}$ the intersection point $\mathcal R_i \cap \mathcal B_j$. 
\smallskip

%\begin{defn} 
Let  $\mathcal A$ be a subset of the nodes of a $(d \times e)$-cage, where $d \geq e$.
We say that $\mathcal A$ is:
\begin{itemize}
\item\emph{triangular} if, with respect to some ordering of the red and blue
lines forming the cage, $\mathcal A$ is of the form  
$\{ p_{ij}\}_{ i + j \leq d}$;
\item\emph{quasi-triangular} if the cardinality of the nodes in $\mathcal A$
residing on a typical blue line ranges from $d$ to $d - e + 1$ 
and takes each intermediate value exactly once;
\item\emph{supra-triangular} if, with
respect to some ordering of the red and blue
lines forming the cage, $\mathcal A$ is of the form 
$\{ p_{ij}\}_{ i + j \leq d + 1}$;
\item\emph{supra-quasi-triangular} if the cardinality of the set of nodes in $\mathcal A$ 
residing on a typical blue line ranges from $d$ to $d - e + 2$, the value 
$d$ is achieved on two lines, and each value from $d - 1$ to $d - e + 2$ is achieved
exactly once. 
\end{itemize}
%\end{defn}
%
\begin{figure}[ht]
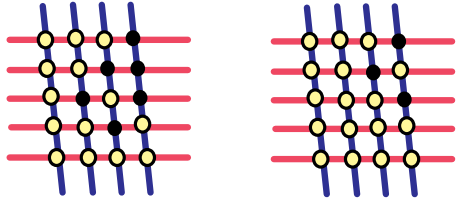

\centerline{\BoxedEPSF{above.triangular scaled 600}} 
\bigskip
\caption{A $(5 \times 4)$-cage with a quasi-triangular and 
an supra-quasi-triangular sets of nodes.}
\end{figure}
We now formulate  our generalization of the $(3\times 3)$-cage
theorem. Later, we will show how this result can be also derived 
with the aid of  Bacharach duality (Theorem 4.1).  
\begin{thm} {\bf(Cage Theorem for Plane Curves).}
 
%\begin{enumerate}
%\item 

{\rm 1.} If a  curve in $\P^2$ of degree $d$ passes through an
supra-quasi-triangular set $\mathcal A$ of nodes of a 
$(d\times e)$-cage with $d \geq e$, then it passes through all the nodes of the cage.

%\item 
{\rm 2.} No curve of degree less than $e$ can pass through a quasi-triangular 
set of nodes of a $(d\times e)$-cage when $d \geq e$.
%\end{enumerate}
\end{thm}
{\bf Remark.} Not every subset of nodes with the  cardinality 
of an supra-triangular set $\mathcal A$ will have the property claimed by 
Theorem 3.1. For example, for a $(4 \times 4)$-cage, $|\mathcal A| = 13$. 
If $\mathcal Y$  is the complement of the set 
 $\{p_{42}, p_{43}, p_{44}\}$ in the nodal set of a $(4 \times 4)$-cage, then 
not every curve $\mathcal C$ of degree four that contains $\mathcal Y$  passes through 
the three nodes $p_{42}, p_{43}$ and $p_{44}$: just take $\mathcal C$ to be
the union of the lines $\mathcal B_1$, $\mathcal B_2$, and $\mathcal B_3$
with a line $\mathcal L$ containing $p_{41}$ but missing $p_{42}, p_{43}$, and $p_{44}$.
Therefore, the combinatorial structure of an supra-quasi-triangular set 
is important for the truth of Theorem 3.1. This example indicates the 
complexity of generalizing the $(4 \times 4)$-cage theorem to sets of 
nodes obtained by intersecting two generic curves of degree four.
\qed 
\smallskip

{\it Proof.} \quad We  prove Theorem 3.1 for cages in $\A^2$.
The argument for cages in $\P^2$ is similar. 
First, for statement (1) we order the blue lines so that the first blue line 
contains $d$ nodes of $\mathcal A$, the second blue line also contains $d$ nodes, the 
third blue line contains $d - 1$ nodes, the fourth $d - 2$ nodes, and so on, 
to the last line, which carries $d - e + 2$ nodes (as depicted in Figure 4, the right-hand diagram). 

Let $R_i$ be a linear polynomial in $x$ and $ y$ 
having the line $\mathcal R_i$ as its zero set,  let $B_i$ be a linear
polynomial with the line $\mathcal B_i$ as its zero set, and let 
$R = \prod_{i=1}^d R_i$ and $B = \prod_{j=1}^e B_j$. 
Then any linear combination $S_{\lambda, \mu} := \lambda R + \mu B$ 
vanishes at each node of the cage. Let $\mathcal S_{[\lambda : \mu]}$ be 
the curve of degree $d$ defined by the equation $S_{\lambda,\, \mu} = 0$. 
This curve depends on the choice of the polynomials 
$R$ and $B$ representing the cage. However, the family $\{\mathcal S_{[\lambda : \mu]}\}$
of such curves is determined completely by the cage. 

Let $P$ be any polynomial of degree at most $d$ that vanishes at the 
nodes from $\mathcal A$. We aim to show that,  
for appropriate numbers $\lambda$ and $\mu$ and some polynomial $Q$,  $P$ is of the form 
$\lambda R + \mu B\cdot Q$,
and therefore must vanish at every node. 

The argument is reductive in nature. Compare the restrictions of $P$ and of 
$S_{\lambda, \mu}$ to the blue line $\mathcal B_1$. Since both restrictions 
are just polynomials of  degree at most $d$ in a single variable (say, the parameter 
$m$ in a parametric description of $\mathcal B_1$) that share a  set of $d$ roots---namely,
the nodes $p_{11}, p_{12}, ... , p_{1d}$---they must be proportional. 
By choosing an appropriate $\lambda = \lambda_\star$, we ensure 
that $P - S_{\lambda_\star,\, \mu}$ is identically zero along  the 
line $\mathcal R_1$. Since $B_1$ is an irreducible polynomial, 
$P - S_{\lambda_\star,\, \mu}$ must be divisible by $B_1$. This statement 
follows from Hilbert's Nullstellensatz (see  [{\bf 10}, Theorem 1.3A]). However, a more elementary 
argument also applies: just pick new local coordinates $(x', y')$ so that $x'$ 
is given by the linear polynomial $B_1$. In the new coordinates, the 
polynomial $P - S_{\lambda_\star,\, \mu}$ transforms into a polynomial $T(x', y')$ in
$x'$ and $ y'$. Since $T(0, y') = 0$ identically in $y'$, all the monomials making up $T$ must 
be divisible by $x' = B_1$. 
Therefore, $P = S_{\lambda_\star,\, \mu} + B_1\cdot P_1$, where $P_1$ is a polynomial 
of degree at most $d - 1$.  Because all the nodes of the set $\mathcal A$ are distinct,
 $P_1$ must vanish at the $d$ nodes of $A \cap \mathcal B_2$ (indeed,
$B_1$ is nonzero on $A \cap \mathcal B_2$). As a consequence, the restriction 
of $P_1$ to the line $\mathcal B_2$ vanishes. Hence, $P_1$ is 
divisible by $B_2$. The original $P$ acquires the form 
$P = S_{\lambda_\star,\, \mu} + (B_1 B_2) \cdot P_2$, where $P_2$ 
is a polynomial of degree at most $d - 2$.  This algorithm can be repeated again 
and again until we reach the conclusion that 
$P = S_{\lambda_\star,\, \mu} + (B_1 B_2 ... B_e) \cdot P_e$, where $P_e$ 
is of degree $d - e$. Hence $P$ must be of the form $\lambda_\star R +  B \cdot Q$. 
In particular, when $e = d$,  we obtain  
$P = S_{\lambda_\star,\, \mu_\star}$ for an appropriate $\mu_\star$. 
\smallskip

The second statement of the cage theorem has an even 
simpler proof. Let $Q$ be a polynomial of 
degree less than $ e$ that vanishes on a quasi-triangular set $T$ of nodes. 
Then, by an argument similar to the one just presented, its restrictions 
to any  blue line must be identical zero: for an appropriate $P_m$, 
the number of roots on $\mathcal B_m$ exceeds its degree. Thus, $Q$ must 
be divisible by $B_1 B_2 ... B_e$, which is only possible when $Q$ is the zero 
polynomial.   \qed
\bigskip
 
Recall that curves of degree $d$ in $\P^2$ are parameterized
by points of the projective space $\P_\ast^{(d^2 + 3d)/2}$ of
dimension $(d^2 + 3d)/2$. The requirement that a curve
of degree $d$ pass through a given point in $\P^2$ imposes a 
homogeneous linear constraint on the homogeneous coordinates 
 in $\P_\ast^{(d^2 + 3d)/2}$. 
Because this  constraint is homogeneous, any system
of $(d^2 + 3d)/2$ such equations has a nontrivial solution. In
other words, for any given $(d^2 + 3d)/2$ points in
$\P^2$ there exists a curve of
degree $d$ that contains them.  

On the other hand, the cardinality of
an supra-quasi-triangular set $\mathcal A$ of nodes  in a $(d\times d)$-cage
is $$d^2 - [(d - 1)(d - 2)/2] = [(d^2 + 3d)/2] - 1,$$
one less than the dimension of the space of degree $d$ curves!
No wonder that the curves that pass through an
supra-quasi-triangular set $\mathcal A$ of a $(d\times d)$-cage form
(according to Theorem 3.2)  a one-parameter family, a line 
$\P_\ast^1$ in $\P_\ast^{(d^2 + 3d)/2}$. The proof of the 
cage theorem tells us that all the $[(d^2 + 3d)/2] - 1$
nodes from $\mathcal A$ impose \emph{independent} linear
restrictions: the nodes of $\mathcal A$
are \emph{generically} located. Evidently, the rest
of the nodes have extremely special locations. The
requirement that a  curve of degree $d$ go through the nodes of
an supra-quasi-triangular set leaves only one degree of freedom.
Furthermore, since any such curve is of the form
$\mathcal{S}_{[\lambda : \mu]}$, the vanishing of the polynomial at
any point in $\P^2$ distinct from the nodes of the cage
fixes the proportionality
class $[\lambda : \mu ]$ and, with it, the corresponding curve in the family 
$\{\mathcal{S}_{[\lambda : \mu]}\}$. 
In particular, when $e = d$, any curve of degree less than $d$ containing an 
supra-quasi-triangular set $\mathcal A$ of nodes and a point on 
a red (blue) line distinct from the nodes is the union 
of red (blue) lines that form the cage.\smallskip

What are some implications of Theorem 3.1?
We have seen that any polynomial of degree at most $d$ that vanishes on an 
supra-quasi-triangular nodal
set $\mathcal A$ of a $(d \times e)$-cage is of the form $\lambda R + B \cdot Q$,
where $\mathrm{deg} \,Q \leq d - e$. Therefore,
for any set $\mathcal F$ of $[(d - e)^2 + 3(d - e)]/2$ points located on the red 
lines of the cage, there is a polynomial $Q$ that vanishes on $\mathcal F$. As a result, the 
zero set of $\lambda R + B\cdot Q$ contains all the nodes of the cage 
together with the set $\mathcal F$.
For a generic set $\mathcal F$ such curves form a one-parameter family. 
This leads to Corollaries 3.1 and 3.2. \smallskip

We say that a finite set $\mathcal B$ is a \emph{minimally redundant set} for 
polynomials of degree $d$ if for some element $b$ of $\mathcal B$
vanishing on the set $\mathcal B \setminus \{b\}$ imposes $|\mathcal B| -1$ independent 
linear constraints  on the coefficients of such polynomials and if any polynomial of degree 
$d$ that vanishes on $\mathcal B \setminus \{b\}$ automatically vanishes on $\mathcal B$. 

\begin{cor} Vanishing requirements at the points of an supra-quasi-triangular 
set $\mathcal A$ of a $(d \times e)$-cage with $e \leq d$ impose
independent conditions on  polynomials of degrees at least $d$. 
Further, adding any new node to $\mathcal A$ produces a 
minimally redundant set $\mathcal B$ that fails to impose independent restrictions 
on polynomials of degree $d$. 
\end{cor}

\begin{cor} If $\mathcal A$ is  an supra-quasi-triangular nodal set  of a $(d \times e)$-cage $K$
and  $\mathcal F$ is a set 
$[(d - e)^2 + 3(d - e)]/2$ points located on the red lines of $K$, then there
exists a one-parameter family of curves of degree $d$ that contain $\mathcal F$ and $\mathcal A$. 
Such  curves  pass  through the rest of the nodes of $K$. 

For a $(d \times d)$-cage $K$, an supra-quasi-triangular nodal set  $\mathcal A$ of $K$, and a point
$p$ different  from the nodes of $K$ there
exists a unique curve $\mathcal C$ of degree $d$ 
that passes through $p$ and  $\mathcal A$. 
The curve  $\mathcal C$ contains all the nodes of $K$. 
\end{cor}

Since at each node $p_{ij}$ of a cage $K$ the red and blue lines $\mathcal R_i$ 
and $\mathcal B_j$ are transversal and since no other lines of the cage pass 
through $p_{ij}$, one can show that any curve  
$\mathcal{S}_{[\lambda : \mu]}$ is \emph{nonsingular} 
at the nodes. It suffices to check that for any $\lambda$ and $\mu$
the gradient  of $S_{\lambda , \mu}$ is nonzero at  $p_{ij}$.
The argument is based on the product rule for gradients. 
The gradient $\nabla S_{\lambda , \mu}$ of $S_{\lambda , \mu}$ at $p_{ij}$ 
is given by 
\begin{equation} 
\nabla S_{\lambda , \mu}(p_{ij}) = 
\lambda r_{ij}\cdot \nabla R_i(p_{ij}) + 
\mu b_{ij} \cdot \nabla B_j(p_{ij}), \nonumber \qquad \qquad \qquad  (1)
\end{equation}
where $r_{ij} = \prod_{k\neq i} R_k(p_{ij})$ and 
$b_{ij} = \prod_{k\neq j} B_k(p_{ij})$ are nonzero. 
Thus the nontrivial linear combination of independent vectors
$\nabla R_i(p_{ij})$ and $\nabla B_j(p_{ij})$ is a nonzero vector.
This implies that each curve  $\mathcal{S}_{[\lambda : \mu]}$ has a 
well-defined tangent line at each of the nodes. Further, as (1) 
testifies, in the cage family there is a \emph{unique} curve tangent to a given 
line $\tau_{ij}^\ast$ in $\P^2_\ast$ that pass through a given node $p_{ij}$. 
Indeed, one can reconstruct $[\lambda:  \mu]$ from the three vectors 
$r_{ij}\cdot \nabla R_1(p_{ij})$,\,  $b_{ij} \cdot \nabla B_1(p_{ij})$, and 
$\nabla S_{\lambda, \mu}(p_{ij})$.    
\smallskip

In this way we can define a  directional line
field $\tau (\lambda , \mu )$ at the nodes of the cage, the one generated 
by the tangent lines to the curve $\mathcal S_{[\lambda : \mu]}$.  
Each direction, say $\tau_{11} (\lambda , \mu )$, 
determines the rest, so that the  field  $\tau (\lambda , \mu )$ 
is ``rigid." This helps to establish the following result:

\begin{thm}  If $K$ is a $(d \times d)$-cage  in $\P^2$ and  
$\tau_{ij}$ in $\P^2_\ast$ is a direction at one of its nodes $p_{ij}$, then there exists a unique curve $\mathcal S_{[\lambda : \mu]}$ of degree  $d$ 
that passes through all the nodes of $K$ and has $\tau_{ij}$ 
as its tangential direction at the  point $p_{ij}$.  

Moreover, the same conclusion holds for any curve of degree $d$  passing 
through an supra-quasi-triangular nodal set $\mathcal A$ of $K$ and having $\tau_{ij}$ as the 
tangential direction at a given point $p_{ij}$ of  $\mathcal A$. 
\end{thm}
\smallskip

We are now in a position to derive a nice generalization (Theorem 3.3) of the
Pascal theorem. It is based simply on Theorem 3.1 and counting dimensions.
\smallskip

The variety of $(d \times d)$-cages is
$4d$-dimensional:  each such cage is determined by a
configuration of $d$ red and $d$ blue points in $\P_\ast^2$. 
In fact, a $(d \times d)$-cage $K$ can be reconstructed from
certain sets consisting of $2d$ nodes. Not any $2d$ nodes will do, 
for in order to determine the cage, they must satisfy certain
combinatorial conditions. One such set  $\mathcal D$
consists of the nodes $p_{ij}$ subject to the constraints
$i + j = d \pm 1$ together with two
``corners" $p_{d1}$ and $p_{1d}$ (see Figure 5).
\begin{figure}[ht]
\centerline{\BoxedEPSF{diagonal scaled 500}}
\bigskip
\caption{ }
\end{figure}
Note that $\mathcal D$ is contained in   
$\mathcal A = \{p_{i j}\}_{i + j \leq d + 1}$, which is an supra-triangular set. 
The complementary set $\mathcal A \setminus \mathcal D$ consists of 
$$\{[(d^2 + 3d)/2] - 1\} - 2d  = [d(d - 1)/2] - 1$$ elements. 
In fact, any subset $\mathcal D$ of $\mathcal A$ with the property  that any line from 
the cage hits it at exactly two points would be a good choice.

Assume for the moment that all $2d$ nodes from the set
$\mathcal D$ are located on an irreducible curve $\mathcal Q$ of
degree $u$ less than $d$. We wish to
determine when the nodes from the set $\mathcal A \setminus \mathcal D$
are located on a curve $\mathcal Q_\star$ of complementary degree $v = d - u$. 
Since there are $[d(d - 1)/2] - 1$ elements in $\mathcal A \setminus \mathcal D$, 
the existence of such a curve is guaranteed by
the inequalities  $$[(v^2 + 3v)/2]  \geq [d(d - 1)/2] - 1$$
and $1 \leq v \leq d - 1$.
Remarkably, when $d$ is at least three, $u = 2$ and $v = d - 2$ satisfy these conditions,
$$[(d - 2)^2 + 3(d - 2)]/2 = [d(d - 1)/2] - 1,$$
whereas any choice of $u$ and $v$ with $u > 2$ fails to obey them! Therefore, the points of $\mathcal A \setminus \mathcal D$ 
are located on a curve $\mathcal Q_\star$ of degree $d - 2$.  
\smallskip

Consider the reducible curve $\mathcal C = \mathcal Q \cup \mathcal Q_\star$ of
degree $d$ with an irreducible quadratic component $\mathcal Q$. The nodes
of $\mathcal D$ are located on $\mathcal Q$, and the nodes from
$\mathcal A \setminus \mathcal D$ are on $\mathcal Q_\star$. Thus, all the nodes
from $\mathcal A$ are located on
$\mathcal C$. By Theorem 3.1, all the nodes of the
cage $K$ must lie on $\mathcal C$ as well. However, the nodes of $K$ that are not
in $\mathcal D$ cannot belong to $\mathcal Q$, because this would contradict
Bezout's theorem (a line would hit the irreducible quadratic curve
$Q$ at three points). It follows that  the
unaccounted for nodes must be located on $\mathcal Q_\star$.

We arrive in this way at our main geometric result---a generalized Pascal's theorem. 
The $2d$-gram that depicts it  is still a bit mysterious, so we retain Pascal's description.

\begin{thm} {\bf (The Mystic 2d-Gram)}. 
Let $\mathcal D$ be a polygon with $2d$ sides, colored in two alternating colors 
and inscribed in a quadratic curve $\mathcal Q$. If $K$ is the $(d \times d)$-cage generated by
$\mathcal D$,\footnote{The vertices of $\mathcal D$ are among the nodes of $K$.}
then all $(d^2 - 2d)$ new nodes of the cage lie on a
plane curve $\mathcal Q_\star$ of degree at most $d - 2$.
For a generic polygon $\mathcal D$, the curve $\mathcal Q_\star$ is unique.
\end{thm}
Figures 6 and 7 show the mystic octagram and decagram, respectively, that 
generate  quadratic and  cubic curves $\mathcal Q_\star$.
\begin{figure}[ht]
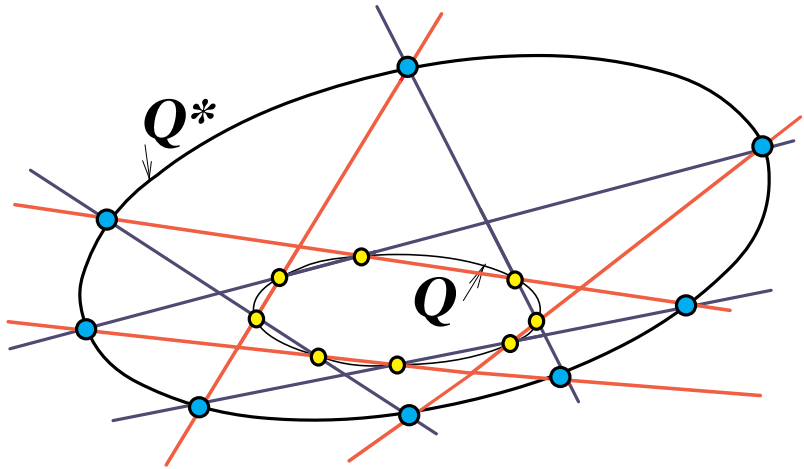

\centerline{\BoxedEPSF{octagon scaled 600}}
\bigskip
\caption{Mystic octagram.}
\end{figure}
\bigskip  
\begin{figure}[ht]
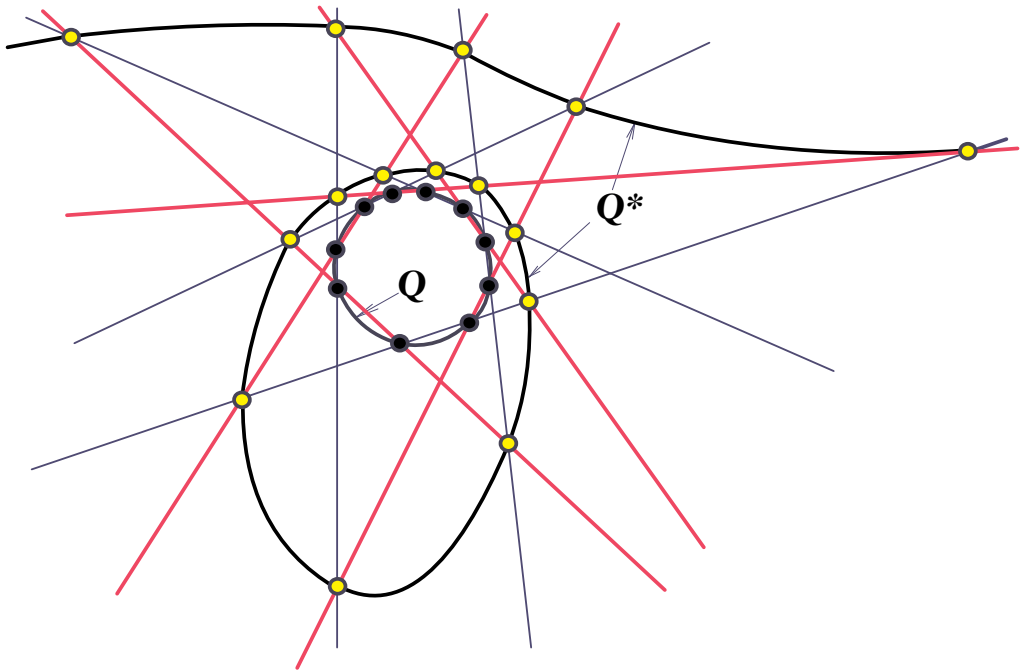

\centerline{\BoxedEPSF{DecagonPascal scaled 700}}
\bigskip
\caption{Mystic decagram.}
\end{figure}
We remark that the $2d$-gon $\mathcal D$ in Theorem 3.3 can be a
union of two or more $2k$-gons
($1 < k < d$), each of which is bicolored in an
alternating manner. Also, as two vertices
merge, some of the sides of the $2d$-gon can become
\emph{tangent} to the quadratic curve.
In particular, if all the red lines become tangent, one obtains
a degenerate cage formed by a red \emph{circumscribed} $d$-gon
and a blue \emph{inscribed} one, the vertices of the blue polygon
being the points of tangency for the red one.
Extending  Theorem 3.3 by continuity, the rest of the nodes of the
resulting cage must lie on a curve of degree $d - 2$.
\smallskip

Contemplating Theorem 3.3, one might wonder:
Which curves $\mathcal Q_\star$ of degree $d - 2$ can
be produced via the mystic $2d$-gon construction from
a $2d$-gon inscribed in
a (given) quadratic $\mathcal Q$? Clearly, for large $d$ such curves
$\mathcal Q_\star$ will be exceptional. However,
for a few small $d$ we might have a chance to
manufacture almost any plane curve of degree $d - 2$
as a $\mathcal Q_\star$. A  count of dimensions  provides a
crude indication to the possible degrees $d$ for which this might be true. 
For $d = 3$ the space of
hexagons inscribed in a particular $\mathcal Q$ is
six-dimensional, and the space of
lines in $\A^2$ is two-dimensional. For $d = 4$ the space of
octagons inscribed in $\mathcal Q$  is eight-dimensional,
while the space of quadratic curves is
five-dimensional. For $d = 5$ the space of  decagons
inscribed in in $\mathcal Q$  is
ten-dimensional, and the space of cubics is
nine-dimensional. So far, so good! But already for $d = 6$
the dimension of the inscribed 12-gons is twelve, whereas the
dimension of quartics is fourteen. Accordingly, not every
quartic can be a $\mathcal Q_\star$ for a \emph{fixed} $\mathcal Q$.
However, if we allow ourselves to vary the quadratic
curve $\mathcal Q$ as well, we gain five extra degrees of freedom.
This  takes us  through the next
case $d = 7$: the space of quintics is 19-dimensional, and
$19 = 14 + 5$. Once  $d$ exceeds seven, the
$\mathcal Q_\star$s  form a subvariety in the space of curves
of degree $d - 2$. 
We conjecture that any plane algebraic curve of 
degree seven or less can be obtained from some $2d$-gon inscribed in 
some quadratic curve via the corresponding $2d$-gram.  
\smallskip 

The next corollary reflects an interesting duality in the universe of
quadratic curves marked with eight ordered points (see Figure 6). 

\begin{cor} {\bf (Mystic Octagram Duality).}
Let $\mathcal D$ be a bicolored octagon inscribed in a quadratic curve $\mathcal Q$.
If $K$ is the $(4 \times 4)$-cage generated by $\mathcal D$,
then  the eight new nodes of the cage lie on a quadratic curve
$\mathcal Q_\star$ and give rise to a new bicolored octagon.
For a generic inscribed octagon $\mathcal D$, the curve $\mathcal Q_\star$
is unique. 
\end{cor}

In particular, the duality tells us that any quadratic
curve $\mathcal Q$ can be generated from another quadratic curve
$\mathcal Q_\star$ via the mystic octagram construction.\smallskip

The octagon in the corollary can be a union of two (bicolored) quadrilaterals. It  can also
degenerate into a union of two (monocolored) quadrilaterals, one of which is tangent to the
quadratic curve and the other of which has its vertices at the points of tangency.
\bigskip

Curiously, quadratic curves play a special role in the mystic
$2d$-gram. As our count of dimensions shows,
\emph{in general}, nothing can be claimed about
 $2d$-gons inscribed in curves of degree $u$ when 
$2 < u < d$. No valid Pascal's theorems inhabit that range!
Nevertheless, the situation is better than one might think, if 
one is  willing to abandon the inscribed $2d$-gons in favor 
of more elaborate inscribed \emph{subcages}. But that is  a subject for another day.
\smallskip 

The following proposition is a well-known consequence of Bezout's theorem 
(see [{\bf 11}, Proposition 3.14]):
\begin{thm} If two  curves $\mathcal D$ and $\mathcal E$ in $\P^2$  of degree $d$ 
intersect in exactly $d^2$ points and if exactly  $k \cdot d$ of these points lie 
on an irreducible curve $\mathcal Q$ of degree $k$ with $k < d$, 
then the remaining $d(d - k)$ intersection points lie on a curve 
$\mathcal Q_\star$ of degree at most $d - k$.
\end{thm}

We concentrate on a special case of Theorem 3.4---namely, the case 
where $\mathcal D$ and $\mathcal E$ 
are the red and blue lines of  a $(d \times d)$-cage $K$. 
Each time we are able to produce  a curve $\mathcal Q$ as 
in Theorem 3.4, a Pascal-type theorem will be  revealed.
However, the unresolved issue raised by Theorem 3.4  
is how to recognize curves $\mathcal D$ and $\mathcal E$ 
whose intersection set admits an irreducible curve $\mathcal Q$  
of degree $k$ passing through $k \cdot d$ nodes.
We do not have a satisfactory answer even  in the case 
when $\mathcal D$ and $\mathcal E$ are each unions of $d$ lines. 
This seems to be a challenging problem. To appreciate the challenge,
consider a $(d \times d)$-cage in $\A^2$ whose nodes have integral coordinates, 
whose red lines are among the horizontal lines forming the integral 
grid, and whose blue lines are among the vertical lines of this grid.
Such a cage is characterized by two unordered $d$-tuples of integers 
$\{n_i\}_{1 \leq i \leq d}$ and $\{m_j\}_{1 \leq j \leq d}$,  the 
$x$- and $y$-intercepts of the blue  and red lines, respectively. 
\smallskip
\begin{figure}[ht]
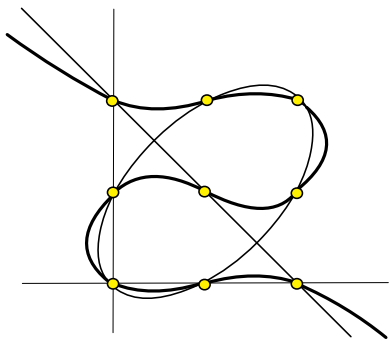

\centerline{\BoxedEPSF{degeneration1 scaled 300}}
\bigskip
\caption{A cubic degenerating into a line and a conic.}
\end{figure}
\begin{figure}[ht]
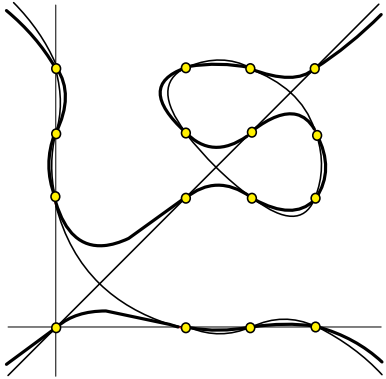

\centerline{\BoxedEPSF{degeneration2 scaled 300}}
\bigskip
\caption{A quartic degenerating into a line and a cubic.}
\end{figure} 
%
%\begin{figure}[ht]
%\centerline{\BoxedEPSF{degeneration3 scaled 300}}
%\bigskip
%\caption{A quartic degenerating into two lines and a quadratic.}
%\end{figure} 
%

{\bf Problem.} Describe in terms of $\{n_i\}$ and $\{m_j\}$ the $(d \times d)$-cages
that admit \emph{reducible} curves of degree $d$ passing through all their nodes 
other than the $\mathcal D$ and $\mathcal E$ forming the cage.
Specifically, we ask: Which of these cages admit  irreducible curves $\mathcal Q$ 
of degree $k$ $( < d)$ that contain $k \cdot d$ nodes?
\smallskip

One class of cages with the desired property is 
easy to produce. These are the cages with collinear ``diagonal"
(with respect to some ordering) nodes $p_{ii}$. 
Figures 8 and 9 display them, together with  the Pascal-type 
diagrams that they support. 
The figures also document how an irreducible curve attached to 
the nodes of such a cage can degenerate into a reducible one.
\smallskip

Theorem 3.4 coupled with the cage theorem implies yet another generalization
of Pascal's theorem:

\begin{cor} The nondiagonal nodes of a $(d \times d)$-cage lie 
on a curve $\mathcal Q$ of degree $d - 1$ if and only if the 
$d$ diagonal nodes of the cage are collinear.  \qed
\end{cor}

%%%%%%%%%%%

\section{cages and modernity}

We close this article by incorporating the elementary considerations of previous sections into  
the language of classical and contemporary algebraic geometry. 
Denote by $\mathcal V_k$ the vector space of homogeneous  polynomials of degree $k$
in $x, y$, and $z$. Its dimension is $(k + 1)(k + 2)/2$. 
Let  $\mathcal X$ be an algebraic subset of  $\P^2$, that is, the solution set for a system of  equations in the variables $x, y$, and $z$ defined by homogeneous polynomials.  We denote by $I_{\mathcal X}$ the ideal of the polynomial ring $\C[x, y, z]$ generated by  polynomials that vanish on $\mathcal X$. 
The vector subspace $I_{\mathcal X,\, k}$ of $I_{\mathcal X}$  comprises its 
homogeneous polynomials of degree $k$. 

The \emph{Hilbert function}  $h_{\mathcal X}:  \Z_+ \rightarrow \Z_+$ of an algebraic set $\mathcal X$ 
is defined by the formula 
$h_{\mathcal X}(k) =  \mathrm {dim} (\mathcal V_k) - \mathrm {dim} (I_{\mathcal X,\, k})$,
where $``\mathrm {dim}"$ signifies the dimension of a vector space.
For sufficiently large $k$,\; $h_{\mathcal X}(k)$ behaves like a polynomial in $k$ of degree 
 $\mathrm {dim}(\mathcal X)$, the complex ``topological" dimension. 
In fact, if $\mathcal X$ is a zero-dimensional set, then for all sufficiently large $k$,
$h_{\mathcal X}(k) = |\mathcal X|$, the cardinality of $\mathcal X$. 
When $\mathcal X$ is a curve,  $h_{\mathcal X}(k)$ is asymptotically a linear polynomial of the form 
$\mathrm{deg}(\mathcal X)\cdot k + c$, where $c$ is a constant (see [{\bf 10}, Theorem 7.5]). \smallskip

The following important theorem describes an intricate duality relation between the Hilbert function
of a complete intersection (roughly speaking, complete intersections are algebraic sets whose codimension in the ambient projective space is equal to the number of polynomial equations that define them) and  the Hilbert functions of its complementary parts {\bf  [1]} (see  
[{\bf 5}, Theorem CB5] and [{\bf 4}] for  significant generalizations):

\begin{thm} {\bf (Bacharach).} Let $\mathcal D$ and $\mathcal E$ be projective plane 
curves of degrees $d$ and $e$, respectively, and let the finite set
$\mathcal X = \mathcal D \cap \mathcal E$ be a complete intersection. 
Assume that $\mathcal X$
is the disjoint union of two subsets $ \mathcal X_1$ and $\mathcal X_2$. Then for any $k$ satisfying 
$k \leq d + e - 3$
\begin{eqnarray}
h_{\mathcal X}(k) - h_{\mathcal X_1}(k) = |\mathcal X_2| - h_{\mathcal X_2}((d + e - 3) - k).
\qquad \qquad \qquad (2) \nonumber
\end{eqnarray}
\end{thm}
The left-hand side of  formula (2)  can be viewed as 
$\mathrm {dim}(I_{\mathcal X_1,\, k} / I_{\mathcal X,\, k})$, while 
the right-hand side measures the failure of $\mathcal X_2$ to impose independent conditions on the polynomials of the ``complementary" degree $(d + e - 3) - k$. \smallskip
Substituting $d$ for $k$ in formula (2), we conclude that 
$h_{\mathcal X}(d) = h_{\mathcal X_1}(d)$ 
if and only if  $|\mathcal X_2| = h_{\mathcal X_2}(e - 3)$. In other words, the following two properties are equivalent: 1) any curve of degree $d$ that contains $\mathcal X_1$ must contain $\mathcal X$ as well, and 2)  $|\mathcal X_2| = h_{\mathcal X_2}(e - 3)$. For instance, with $d = e = 3$, we get the 
Chasles and $(3\times 3)$-cage theorems; if $|\mathcal X_2| = 1$, then $h_{\mathcal X_2}(3 - 3) = 1$ 
(a polynomial of degree zero that vanishes at the singleton  $\mathcal X_2$ must be zero),
which by the Bacharach theorem implies that $h_{\mathcal X}(3) = h_{\mathcal X_1}(3)$. In other words, $I_{\mathcal X, 3} = I_{\mathcal X_1, 3}$, so any cubic polynomial that 
vanishes on $\mathcal X_1$ of cardinality eight vanishes on $\mathcal X$  of cardinality nine.
\smallskip

We now see how  Theorem 3.1 can be derived from Theorem 4.1.  To simplify the combinatorics, we 
consider only $(d \times d)$-cages and only triangular and supra-triangular sets of nodes. 
We work by induction on $d$. 
Assume that the cage theorem (which consists of two claims) is true whenever 
$d < d_\star$. In particular, this inductive assumption implies that, when $d < d_\star$,
any polynomial of degree less than $d$ that vanishes at the points of a triangular set in 
a $(d \times d)$-cage is identically zero.

Let $\mathcal A$ denote an supra-triangular set of nodes 
in a $(d_\star \times d_\star)$-cage, and let $\mathcal B$ be its complementary 
set. The first statement of the cage theorem can be expressed as 
 $\mathcal I_{\mathcal A, d_\star} = \mathcal I_{\mathcal A \cup \mathcal B, d_\star}$ which 
is equivalent to the claim that $h_{\mathcal A}(d_\star) = h_{\mathcal A\cup \mathcal B}(d_\star)$.  
By  Theorem 4.1, this  is equivalent to the property
$h_{\mathcal B}(d_\star - 3) = |\mathcal B|$. Note that $\mathcal B$ is a triangular set 
for an appropriate $(d_\star - 2)\times (d_\star - 2)$-subcage. 
By the inductive assumption, any polynomial of degree $d_\star - 3$ 
that vanishes on $\mathcal B$ is the zero polynomial, that is, 
$\mathcal I_{\mathcal B, d_\star - 3} = \{0\}$.
Hence $h_{\mathcal B}(d_\star - 3) = 
\mathrm{dim}(\mathcal V_{d_\star - 3})$. On the other hand, $\mathrm{dim}(\mathcal V_{d_\star - 3}) = 
(d_\star - 1)(d_\star - 2)/2 = |\mathcal B|$, 
which completes the inductive step for the first statement of the 
cage theorem.  

The inductive treatment of the second statement is similar to the one presented in Theorem 3.1.
Consider a polynomial $P$ of degree $d_\star - 1$ that vanishes on a triangular set $\mathcal T$ of 
a $(d_\star \times d_\star)$-cage. Then it vanishes at the $d_\star$ nodes of the cage 
that lie on the line $\mathcal B_1: B_1 = 0$. As a result, $P$ must be divisible by $B_1$. 
Let $P = B_1\cdot Q$. Then $Q$ is of degree $d_\star - 2$ and must vanish at the nodes of $\mathcal T$ 
that do not belong to $\mathcal B_1$. They form a triangular set of a 
$((d_\star - 1)\times (d_\star - 1))$-cage. Thus, by induction, $Q$ is the zero polynomial, and so must be 
$P$. \qed \bigskip

Although in our investigations we have managed without Hilbert functions, they are invaluable 
in  modern research relevant to the subject of our paper 
(see [{\bf 6}], [{\bf 7}], [{\bf 8}], [{\bf 9}]).   
The interested reader may examine [{\bf 9}, Theorem 3.13] to get a feel 
for the most recent developments in the field. This comprehensive result 
by Geramita, Harita, and Shin is a Fubini-type theorem for the Hilbert function 
of a finite subset $\mathcal X$ of $\P^n$ that  is contained in the union of a family 
of  hypersurfaces $\{\mathcal V_i\}_{1 \leq i \leq s}$ whose  degrees $d_i$  add up to the 
degree of $\mathcal X$. Under a subtle hypothesis that regulates 
the interaction between $\mathcal X$ and the hypersurfaces, a nice 
formula emerges:
$$h_{\mathcal X}(k) =  h_{\mathcal X \cap \mathcal V_1}(k) + 
h_{\mathcal X \cap \mathcal V_2}(k - d_1) + \cdots
  +  h_{\mathcal X \cap \mathcal V_s}(k - (d_1 + \cdots  + d_{s-1})).$$
This formula could be instrumental in generalizing our results to 
multidimensional cages, or even to  $\mathcal X$s that are not necessarily 
complete intersections.  
\bigskip 

{\bf  ACKNOWLEDGMENTS.} I am grateful to the referee, whose suggestions helped to 
improve the quality and scope of the original manuscript.  In particular, his remarks 
helped me to clarify the relation between the cage theorem and its  
generalizations at the very center of more advanced and active research. My  gratitude also
goes to my colleagues Harry Tamvakis and Glen Van Brummelen, who helped me achieve a better exposition.

\bigskip

{\bf GABRIEL KATZ} grew up in Kishinev, Moldova. He received his M.A. and Ph.D. degrees from  Moscow State University. After immigrating to Israel in 1979, he taught at Tel Aviv University and later at Ben Gurion University, where he held a tenured position. In the early '90s, he moved with his family to the United States and settled down in the Boston area. Since then, he has held  a variety of teaching and research positions  at a number of  universities and colleges:  Rutgers, Brandeis, Harvard,  MIT, Clark, Wellesley,  and Bennington, to name a few. Currently, he is a research affiliate at  Brandeis University. 
In the first half of his academic carrier, Dr. Katz  studied  symmetries of  smooth manifolds.  Lately, his work has acquired a more geometrical flavor;  his articles include topics that cover:  the topology of harmonic one-forms,  three-dimensional topology, algebraic geometry of discriminant varieties,  and more recently,  the Morse Theory on manifolds with boundary.  Dr. Katz also does research in the field of Mathematics Education.

\end{document}